\documentclass[a4paper,11pt]{article} 

\usepackage{amsmath, amssymb, comment}
\usepackage{shadow, epsfig, color}
\usepackage{graphicx}
\usepackage{url}
\usepackage{mathrsfs}

\newcommand{\N}{\mathbb N}

\newcommand{\R}{\mathbb R}

\newcommand{\e}{\varepsilon}
\newcommand{\1}{\mathbf 1}
\newcommand{\qed}{\ \hfill \fbox{} \bigskip}
\newcommand{\proof}{{\it Proof.} }

\newcommand{\E}{\mathcal E}
\newcommand{\D}{\mathcal D}
\newcommand{\F}{\mathcal F}

\numberwithin{equation}{section}

\newtheorem{thm}{Theorem}[section]
\newtheorem{defn}[thm]{Definition}
\newtheorem{lem}[thm]{Lemma}
\newtheorem{rem}[thm]{Remark}
\newtheorem{cor}[thm]{Corollary}

\newtheorem{prop}[thm]{Proposition}
\newtheorem{asmp}[thm]{Assumption}

\makeatletter
\newcommand{\subjclass}[2][2010]{%
  \let\@oldtitle\@title%
  \gdef\@title{\@oldtitle\footnotetext{#1 \emph{Mathematics subject classification.} #2}}%
}
\newcommand{\keywords}[1]{%
  \let\@@oldtitle\@title%
  \gdef\@title{\@@oldtitle\footnotetext{\emph{Key words and phrases.} #1.}}%
}
\makeatother

\title{Convergence of continuous stochastic processes on compact metric spaces converging in the Lipschitz distance}

\author{Kohei Suzuki \vspace{2mm} \\ {\it \small Department of Mathematics, Faculty of Science} \\ {\it \small Kyoto University} 
   \\ {\it \small Kyoto, 606-8502, Japan}}


\begin{document}
\date{}
\maketitle

\begin{abstract}
We introduce a new distance, {\it a Lipschitz--Prokhorov distance} $d_{LP}$, on the set $\mathcal {PM}$ of isomorphism classes of pairs $(X, P)$ where $X$ is a compact metric space and $P$ is the law of a continuous stochastic process on $X$.
We show that $(\mathcal {PM}, d_{LP})$ is a complete metric space. 
For Markov processes on Riemannian manifolds, we study relative compactness and convergence.
\end{abstract}

\section{Introduction}
The motivation of this paper is to study a convergence of continuous stochastic processes on varying compact metric spaces.
If state spaces are fixed, we have the weak convergence as a standard notion of convergences of stochastic processes.
When state spaces are not fixed, but embedded into one common space, it is possible to consider the weak convergence.
There have been many studies in such embedded situations: for example, 
approximations of diffusion processes on $\mathbb R^d$ by discrete Markov chains on $(1/n)\mathbb Z^d$ in \cite{SV79}, \cite{SZ97} and \cite{BC08}; approximations of jump processes on proper metric spaces by Markov chains on discrete graphs in \cite{BKU10} and \cite{CKK13}, and on ultra-metric spaces in \cite{S14}; diffusion processes on thin tubes shrinking to graphs in $\mathbb R^d$ in \cite{AK12} and references therein; many studies about scaling limits of random processes on random environments (see, e.g., \cite{K14} and references therein).

In this paper, we consider a convergence of continuous stochastic processes on compact metric spaces converging {\it in the Lipschitz distance.} In this setting, state spaces are not necessarily embedded initially into one common space. By the aid of the Lipschitz convergence, however, we can choose a family of bi-Lipschitz embeddings, which enables us to embed varying state spaces into one common space. After embedding, we can consider the weak convergence of stochastic processes on such the common space. With such ideas, we will first introduce a new distance on the set of pairs of compact metric spaces and continuous stochastic processes,  and show the completeness of this distance. We will second study several topological properties induced by this new distance. 

To be more precise, the main object of this paper is a pair $(X,P)$ where $X$ is a compact metric space and $P$ is the law of a continuous stochastic process on $X$. 
Let $\mathcal {PM}$ be the set of all pairs $(X,P)$ modulo by an isomorphism relation (defined in Section \ref{sec: LP}).  We will define (in Section \ref{sec: LP}) a new distance on $\mathcal {PM}$, which we will call {\it the Lipschitz--Prokhorov distance} $d_{LP}$,  as a kind of mixture of the Lipschitz distance and the Prokhorov distance. The Lipschitz distance is  a distance on the set of isometry classes of metric spaces, which was first introduced by Gromov (see e.g., \cite{Gro99}).

We summarize our results as follows:
\begin{description}
	\item[(A)] $(\mathcal {PM}, d_{LP})$ is a complete metric space (Theorem \ref{thm: LPD} and Theorem \ref{thm: Polish});
	\item[(B)] Relative compactness in $(\mathcal {PM}, d_{LP})$ follows from bounds for sectional curvatures, diameters and volumes of Riemannian manifolds and uniform heat kernel estimates of Markov processes (Theorem \ref{thm: RRC});
	\item[(C)] Sequences in a relatively compact set are convergent if the corresponding Dirichlet forms of Markov processes are Mosco-convergent in the sense of Kuwae--Shioya \cite{KS03} (Theorem \ref{thm: SCF});
	\item[(D)] Examples for 
	\begin{itemize}
		\item Brownian motions on Riemannian manifolds (Section \ref{subsec: BM});
		\item uniformly elliptic diffusions on Riemannian manifolds (Section \ref{subsec: UE}).
	\end{itemize}
\end{description}

Let us explain (A). Let $\mathcal C(X)$ be the set of continuous paths from $[0,T]$ to $X$ equipped with the uniform metric where $T>0$ is a fixed positive real number. A map $f: (X,P) \to (Y,Q)$ is called {\it an $(\e,\delta)$-isomorphism} 
if 
\begin{itemize}
	\item $f:X\to Y$ is an $\e$-isometry (see Section \ref{sec: LP});
	\item The following inequalities hold:
	\begin{equation}\label{ineq: LPD1}
	\begin{split}
& {\Phi_f}_*P(A) \le Q(A^{\delta e^\e})+\delta e^\e, \ Q(A) \le {\Phi_f}_*P(A^{\delta e^\e})+\delta e^\e,
\\
&{\Phi_{f^{-1}}}_*Q(B) \le P(B^{\delta e^\e})+\delta e^\e, \ P(B) \le {\Phi_{f^{-1}}}_*Q(B^{\delta e^\e})+\delta e^\e
\end{split}
	\end{equation}
for any Borel sets $A \subset \mathcal C(Y)$ and $B \subset \mathcal C(X)$ and we mean that  $\Phi_f: \mathcal C(X) \to \mathcal C(Y)$ is defined by $v \mapsto f(v(t))$.
\end{itemize} 
 We define $d_{LP}$ as 
\begin{align} \label{defn-equa: LPD}
	d_{LP}((X,P),(Y,Q))=\inf\{\e+\delta \ge 0: \exists (\e,\delta)\text{-isomorphism} \}.
\end{align}

The inequalities \eqref{ineq: LPD1} indicate how to measure the distance between $P$ and $Q$, which live on different path spaces $\mathcal C(X)$ and $\mathcal C(Y)$: First we push-forward $P$ by $\Phi_f$, and then ${\Phi_f}_*P$ and $Q$ live on the same path space $\mathcal C(Y)$. Second we measure ${\Phi_f}_*P$ and $Q$ by a kind of a modified Prokhorov metric, which involves a space error $e^\e$ due to an $\e$-isometry. Note that, if we replace $e^\e$ to $1$ in \eqref{ineq: LPD1}, then the triangle inequality for the Lipschitz--Prokhorov distance $d_{LP}$ fails.

We explain (B). After we have the complete metric space $(\mathcal {PM}, d_{LP})$, one of the important questions is:
\begin{description}
\item[(Q1)] Which subsets in $(\mathcal {PM}, d_{LP})$ are compact?
\end{description}
We restrict our interest to paris of Riemannian manifolds and Markov processes. 
We introduce a certain subset $\mathcal P_\phi \mathcal R=\mathcal P_\phi \mathcal R(n,K,V,D) \subset \mathcal{PM }$
consisting of pairs $(M,P)$ where $M$ is a Riemannian manifold with bounds for the sectional curvature, diameter and volume, and $P=P^\mu$ be the law of a Markov process with an initial distribution $\mu$ associated with a Dirichlet form whose heat kernel has a uniform bound by a given function $\phi$ (see details in Section \ref{sec: RC}). In Theorem \ref{thm: RRC}, we will show that $\mathcal P_\phi \mathcal R$ is relatively compact in $(\mathcal {PM}, d_{LP})$. 

We explain (C). Let $(M_i,P_i)$ be a sequence in the relatively compact subset $\mathcal P_\phi \mathcal R$. By the relative compactness, we can take a converging subsequence from a sequence $(M_i, P_i) \in \mathcal P_\phi \mathcal R$.
Then a question is: 
 \begin{description}
\item[(Q2)] When does $(M_i, P_i)$ converge in $d_{LP}$ without taking subsequences?
\end{description}
In Theorem \ref{thm: SCF}, we will give a sufficient condition for such convergence in terms of the Mosco-convergence of the corresponding Dirichlet forms in the sense of Kuwae-Shioya \cite{KS03}. The Mosco-convergence was first introduced by Mosco \cite{M67} (see also \cite{M94}). In \cite{KS03}, they generalized the Mosco-convergence to the case of varying state spaces.

Here we refer to some related topics. There have been many studies of convergences of analytical objects related to Markov processes, such as eigenvalues of Laplacians, heat kernels, Dirichlet forms, tensor fields, differentials of Lipschitz functions, on state spaces converging in the (measured) Gromov--Hausdorff sense in \cite{F87, KK94, KK96, KKO97, KMS01, S01, K02, KK02, K06, KS03, KS08, H11, H13a, H13b, H14}. 
In Ogura \cite{O01}, the author dealt with a convergence of stochastic processes on varying state spaces adopting rather a different approach from ours.
He considered a convergence of Brownian motions on Riemannian manifolds under the Gromov--Hausdorff convergence.
Since approximation maps are not necessarily continuous in the case of the Gromov--Hausdorff convergence, the push-forward measures of the laws of continuous stochastic processes do not necessarily live on 
the continuous path space (but live on the Borel measurable path space). This constitutes a great difficulty because there is no standard notion of convergence for probability measures on the Borel measurable path space. In \cite{O01}, the author overcame this difficulty by discretizing the time parameters of processes, which makes the push-forward
measures live on the space of right-continuous paths with left-hand limits. 

The present paper is organized as follows. In Section \ref{sec: LP}, we introduce the Lipschitz--Prokhorov distance $d_{LP}$ and show that $d_{LP}$ is complete.
In Section \ref{sec: RC}, we give a sufficient condition for relative compactness and also give a sufficient condition for sequences in relatively compact sets to be convergent without taking subsequences.
In Section \ref{sec: EX}, we give several examples. In Section \ref{subsec: BM}, we consider the case of Brownian motion on Riemannian manifolds. In Section \ref{subsec: UE}, we consider the case of 
diffusions associated with the uniformly elliptic second-order differential operators.

\section{Lipschitz--Prokhorov distance} \label{sec: LP}

In this section, we introduce a distance $d_{LP}$ on $\mathcal {P}\mathcal M$, called {\it the Lipschitz--Prokhorov distance} and show $(\mathcal {PM}, d_{LP})$ is a complete metric space.
We first recall the Lipschitz distance and the Prokhorov distance briefly. 

Recall the Lipschitz distance. We say that a map $f:X \to Y$ between two metric spaces is {\it an isometry} if $f$ is surjective and distance preserving. Let $\mathcal M$ denote the set of isometry classes of compact metric spaces.
Let $X$ and $Y$ be in $\mathcal M$. For a bi-Lipschitz homeomorphism $f: X \to Y$, {\it the dilation of $f$} is defined to be the smallest Lipschitz constant
of $f$:
$${\rm dil}(f)=\sup_{x\neq y}\frac{d(f(x),f(y))}{d(x,y)}.$$ For $\e\ge 0$, a bi-Lipschitz homeomorphism  $f: X \to Y$ is said to be {\it an $\e$-isometry} if 
$$|\log {\rm dil}(f)|+|\log {\rm dil}(f^{-1})| \le \e.$$
By definition, $0$-isometry is an isometry. 
{\it The Lipschitz distance $d_L(X, Y)$} between $X$ and $Y$ is defined to be the infimum of $\e \ge 0$ such that an $\e$-isometry between $X$ and $Y$ exists:
$$d_L(X,Y)=\inf\{\e \ge 0: \exists f: X\to Y \ \text{$\e$-isometry}\}.$$
If no bi-Lipschitz homeomorphism exists between $X$ and $Y$, we define $d_L(X,Y)=\infty$.
We say that a sequence $X_i$ in $\mathcal M$ {\it Lipschitz converges to }$X$ if $$d_L(X_i,X) \to 0 \quad (i \to \infty).$$
We note that $(\mathcal M, d_L)$ is a complete metric space. This fact may be known, but we do not know any references and, for the readers' convenience, we will prove the completeness of $(\mathcal M, d_L)$ in Proposition \ref{prop: NSC} in Appendix.
Note that $(\mathcal M, d_L)$ is not separable because the Hausdorff dimensions of $X$ and $Y$ must coincide if $d_L(X, Y)<\infty$. See Remark \ref{rem: NS}. We refer the reader to e.g., \cite{Gro99, BBI01} for details of the Lipschitz convergence.

Recall the Prokhorov distance.  For $T>0$, let $\mathcal C_T(X)$ denote the space of continuous maps from $[0,T]$ to a compact metric space $X$ with the uniform metric 
$$d_{\mathcal C}(v, w)= \sup_{t \in [0,T]}d(v(t), w(t)) \quad \text{for} \quad v, w \in \mathcal C_T(X).$$ We fix a constant $T>0$ and we write $\mathcal C(X)$ shortly for $\mathcal C_T(X)$. 
Let $\mathcal P(\mathcal C(X))$ denote the set of probability measures on $\mathcal C(X)$. 
{\it The Prokhorov distance} between two probability measures $P$ and $Q$ on $\mathcal C(X)$ is defined to be 
\begin{align*}
	d_P(P,Q)=\inf\{\delta\ge 0: & P(A) \le Q(A^\delta)+\delta, \ Q(A) \le P(A^\delta)+\delta 
	\\
	&\text{for any Borel set $A\subset \mathcal C(X)$}\},
\end{align*}
where $A^\delta=\{x \in \mathcal C(X): d_\mathcal C(x,A) < \delta\}$.
We know that $(\mathcal P(\mathcal C(X)), d_P)$ is a complete separable metric space (see \cite[\S 6]{B99}).
We refer the reader to e.g., \cite[\S 6]{B99} for details of the Prokhorov distance.

Now we introduce {\it the Lipschitz--Prokhorov distance}. 
For a continuous map $f: X\to Y$, we define $\Phi_f: \mathcal C(X) \to \mathcal C(Y)$ by $$\Phi_f(v)(t)=f(v(t)) \quad (v \in \mathcal C(X), t\in [0,T]).$$  Let $(X,P)$ be a pair of a compact metric space $X$ and a probability measure $P$ on $\mathcal C(X)$. Note that $$\text{$P$ is {\it not} a probability measure {\it on} $X$, but {\it on} $\mathcal C(X)$}.$$ We say that two pairs of $(X, P)$ and $(Y, Q)$ are {\it isomorphic} if there is an isometry $f:X \to Y$ such that the push-forward measure ${\Phi_{f}}_*P$ is equal to $Q$. Note that ${\Phi_{f}}_*P=Q$ implies ${\Phi_{f^{-1}}}_*Q=P$ and thus the isomorphic relation becomes an equivalence relation. 
Let $\mathcal {P}\mathcal M$ denote the set of isomorphism classes of pairs $(X, P)$. 
Let $(X, P)$ and $(Y, Q)$ be in $\mathcal {P}\mathcal M$. Now we introduce a notion of {\it an $(\e, \delta)$-isomorphism}, which is a kind of generalization of $\e$-isometry. A map $f: (X,P) \to (Y,Q)$ is called {\it an $(\e, \delta)$-isomorphism} if the following hold:
\begin{description}
\item[(i)] $f:X\to Y$ is an $\e$-isometry;
\item[(ii)] the following inequalities hold:
\begin{equation}\label{eq: LPD}
	\begin{split}
& {\Phi_f}_*P(A) \le Q(A^{\delta e^\e})+\delta e^\e, \ Q(A) \le {\Phi_f}_*P(A^{\delta e^\e})+\delta e^\e,
\\
&{\Phi_{f^{-1}}}_*Q(B) \le P(B^{\delta e^\e})+\delta e^\e, \ P(B) \le {\Phi_{f^{-1}}}_*Q(B^{\delta e^\e})+\delta e^\e,
\end{split}
	\end{equation}
for any Borel sets $A \subset \mathcal C(Y)$ and $B \subset \mathcal C(X)$.
\end{description}
We now define a distance between $(X, P)$ and $(Y, Q)$ in $\mathcal P\mathcal M$, which is called {\it the Lipschitz--Prokhorov distance}.
\begin{defn} \label{defn: LPD} \normalfont
Let $(X, P)$ and $(Y, Q)$ be in $\mathcal {P}\mathcal M$. {\it The Lipschitz--Prokhorov distance} between $(X, P)$ and $(Y, Q)$ is defined to be the infimum of $\e+\delta\ge 0$ such that an $(\e, \delta)$-isomorphism $f:(X, P) \to (Y, Q)$ exists:
\begin{align*}
	d_{LP}&\bigl((X, P), (Y, Q)\bigr)
=\inf\{\e+\delta \ge 0: \exists f: (X,P)\to (Y,Q) \ \text{$(\e,\delta)$-isomorphism}\}.
\end{align*}
If there is no $(\e, \delta)$-isomorphism between $(X,P)$ and $(Y,Q)$, we define $$d_{LP}\bigl((X,P), (Y,Q)\bigr)=\infty.$$ 
\end{defn}
\begin{rem} \normalfont
If we replace $e^{\e}$ to $1$ in the inequalities \eqref{eq: LPD}, then the triangle inequality fails for $d_{LP}$.
\end{rem}
\begin{rem} \normalfont 
If we start at metric measure spaces, it may seem to be more natural than $(X,P)$ to consider triplets $(\mathcal C(X), d_{\mathcal C}, P)$ where $X \in \mathcal M$ and $P$ is a probability measure on $\mathcal C(X)$.
It is, however, not suitable for our motivation because the Lipschitz convergence of $\mathcal C(X_i)$ does not imply the Lipschitz convergence of $X_i$ in general. 
\end{rem}

It is clear by definition that $d_{LP}$ is well-defined in $\mathcal P\mathcal M$, that is, if $(X, P)$ is isomorphic to $(X', P')$ and $(Y,Q)$ is isomorphic to $(Y', Q')$, then $$d_{LP}\bigl((X, P), (Y, Q)\bigr)=d_{LP}\bigl((X', P'), (Y', Q')\bigr).$$ It is also clear by definition that $d_{LP}((X, P), (X, P))=0$, and $d_{LP}$ is non-negative and symmetric.
To show that $d_{LP}$ is a metric on $\mathcal P\mathcal M$, it is enough to show that $d_{LP}$ satisfies the triangle inequality and that $(X, P)$ and $(Y, Q)$ are isomorphic if $d_{LP}\bigl((X, P), (Y, Q)\bigr)=0$.
Before the proof, we utilize the following lemma:
\begin{lem} \label{lem: isometry}
Let $f:X\to Y$ be an $\e$-isometry. Then $\Phi_f: \mathcal C(X) \to \mathcal C(Y)$ is also an $\e$-isometry with respect to the uniform metric $d_{\mathcal C}$. As a byproduct,  for any $a\ge 0$ and Borel set $A \subset \mathcal C(Y)$,  we have 
$$\Phi_f^{-1}(A^a) \subset \Phi_f^{-1}(A)^{a e^\e} \quad \text{and} \quad \Phi_f^{-1}(A)^a \subset \Phi_f^{-1}(A^{a e^\e}), $$
and, for any Borel set $B \subset \mathcal C(X)$, we have 
$$\Phi_f(B^a) \subset \Phi_f(B)^{a e^\e} \quad \text{and} \quad \Phi_f(B)^a \subset \Phi_f(B^{ae^\e}).$$
\end{lem}
\proof
We first show that $\Phi_f: \mathcal C(X) \to \mathcal C(Y)$ is an $\e$-isometry.
It is clear that $\Phi_f$ is a homeomorphism.
Let $v,w \in \mathcal C(X)$. By the compactness of $[0,T]$ and the continuity of $v,w$ and $f$, there are $t_0 \in [0,T]$ and $s_0 \in [0,T]$ such that 
$$d_{\mathcal C} (v, w)=d\bigl(v(t_0), w(t_0)\bigr) \quad \text{and} \quad  d_{\mathcal C}\bigr(\Phi_f(v), \Phi_f(w)\bigl)=d\Bigl(f\bigl(v(s_0)\bigr), f\bigl(w(s_0)\bigr)\Bigr).$$
Then we have 
\begin{align*}
\frac{d_{\mathcal C}\bigl(\Phi_f(v), \Phi_f(w)\bigr)}{d_{\mathcal C}(v,w)} &= \frac{d\Bigl(f\bigl(v(s_0)\bigr), f\bigl(w(s_0)\bigr)\Bigr)}{d\bigl(v(t_0), w(t_0)\bigr) }
\\
& \le {\rm dil}(f) \frac{d\bigl(v(s_0), w(s_0)\bigr)}{d\bigl(v(t_0), w(t_0)\bigr) }
\\
& \le {\rm dil}(f).
\end{align*}
We also have 
\begin{align*}
\frac{d\bigl(f(x), f(y)\bigr)}{d(x,y)} = \frac{d_{\mathcal C}\bigl(\Phi_f(c_x), \Phi_f(c_y)\bigr)}{d_{\mathcal C}(c_x, c_y) }
 \le {\rm dil}(\Phi_f).
\end{align*}
Here $c_{x}$ denotes the constant path on $x \in X$, that is, $c_{x}(t)=x$ for all $t \in [0,T]$.
Thus we have ${\rm dil}(\Phi_f) = {\rm dil}(f)$. By the same argument, we also have ${\rm dil}(\Phi^{-1}_f) = {\rm dil}(f^{-1})$. These imply that  $\Phi_f$ is an $\e$-isometry.

We second show the inclusions in the statement. It is enough to show one of the inclusions, say, $\Phi_f(B^a) \subset \Phi_f(B)^{a e^\e}$ for $a \ge 0$ and  any Borel sets $B \subset \mathcal C(X)$.
Let $x \in \Phi_f(B^a)$ and $y \in B^a$ such that $ \Phi_f(y)=x$.
Since $\Phi_f$ is an $\e$-isometry, we have 
\begin{align*}
	d_\mathcal C (x, \Phi_f(B)) \le {\rm dil}(\Phi_f)d_\mathcal C (y, B) \le e^\e d_\mathcal C (y, B) \le a e^\e.
\end{align*}
Thus we have $x \in \Phi_f(B)^{a e^\e}$ and finish the proof.
\qed

Now we show that $d_{LP}$ is a metric on $\mathcal P\mathcal M$.
\begin{thm} \label{thm: LPD}
$d_{LP}$ is a metric on $\mathcal {P}\mathcal M$.
\end{thm}
\proof
It is enough to show the following two statements:
\begin{description}
	\item[{(i)}] $d_{LP}$ satisfies the triangle inequality;
	\item[(ii)] $d_{LP}\bigl((X, P), (Y, Q)\bigr) =0$ implies that 
$(X, P)$ and $(Y, Q)$ are isomorphic.
\end{description}

We first show the statement {\bf (i)}. Let $(X, P), (Y, Q)$ and $(Z, R) \in \mathcal {P}\mathcal M$ such that 
there are $(\e_1, \delta_1)$-isomorphism $f_1: X \to Y$ and $(\e_2, \delta_2)$-isomorphism $f_2: Y \to Z$.    
It suffices to show that $f_2\circ f_1:X \to Z$ is an $(\e_1+\e_2, \delta_1+\delta_2)$-isomorphism. In fact, this implies 
$$d_{LP}\bigl((X, P), (Z, R)\bigr)<\e_1+\e_2+\delta_1+\delta_2.$$ 
By taking the infimum of $\e_1+\delta_1$ and $\e_2+\delta_2$, we have the triangle inequality.

Thus we now show that $f_2\circ f_1:X \to Z$ is an $(\e_1+\e_2, \delta_1+\delta_2)$-isomorphism. We know that $f_2\circ f_1$ is an $(\e_1 + \e_2)$-isometry (see e.g., \cite[Theorem 7.2.4]{BBI01} ).
For any Borel set $A \subset \mathcal C(Z)$, we have
\begin{align*}
	{\Phi_{f_2\circ f_1}}_*P(A) &=({\Phi_{f_1}}_*P)(\Phi_{f_2}^{-1}(A))
	\\
	&\le Q\bigl(\Phi_{f_2}^{-1}(A)^{\delta_1e^{\e_1}}\bigr)+\delta_1e^{\e_1}
	\\
	&\le \textcolor{black}{Q\bigl(\Phi_{f_2}^{-1}(A^{\delta_1e^{\e_1+\e_2}})\bigr)+\delta_1e^{\e_1} }
	\\
	&\le R(A^{\delta_1 e^{ \e_1+\e_2}+\delta_2 e^{\e_2}})+\delta_1e^{\e_1}+\delta_2e^{\e_2}
	\\
	& \le R(A^{(\delta_1+\delta_2) e^{ \e_1+\e_2}})+(\delta_1+\delta_2) e^{ \e_1+\e_2}.
\end{align*}
The inequality of the third line follows from Lemma \ref{lem: isometry}. The other directions of \eqref{eq: LPD} can be shown by the same argument. Thus we have that $f_2\circ f_1$ is an $(\e_1+\e_2, \delta_1+\delta_2)$-isomorphism. We have the triangle inequality.

We second show {\bf (ii)}. We show that there is an isometry $\iota: X \to Y$ such that $\Phi_\iota: (\mathcal C(X), P) \to (\mathcal C(Y), Q)$ is a measure-preserving map, that is, for any real-valued uniformly continuous and bounded 
function $u$ on $\mathcal C(Y)$,  we have 
\begin{align} \label{eq: isom}
\int_{\mathcal C(X)}u\circ \Phi_\iota \ dP=\int_{\mathcal C(Y)}u \ dQ.
\end{align}
Let $f_i: X \to Y$ be an $(\e_i, \delta_i)$-isomorphism with $\e_i, \delta_i \to 0$ as $i \to \infty$.
Since the dilations of $\{f_i: i \in \N\}$ are uniformly bounded, by the Ascoli--Arzel\`a theorem, we can take a subsequence from $\{f_i: i \in \N\}$ converging uniformly to a continuous function $\iota$.
Since $\e_i\to 0$ as $i \to \infty$, we have that $\iota$ is an isometry from $X$ to $Y$, and, for any $\e>0$, there is an $i_0$ such that, for any $i_0 \le i$, we have $d(\iota(x), f_i(x))<\e$ for all $x \in X$. See e.g., 
\cite[Theorem 7.2.4]{BBI01} for details. By this fact, we have  
$$d_{\mathcal C}(\Phi_\iota(v),\Phi_{f_i}(v)) \le \e  \quad (\forall i \ge i_0, \forall v \in \mathcal C(X)).$$
By the uniform continuity of $u$, we have 
\begin{align} \label{eq: ISO-1}
\Bigl| \int_{\mathcal C(X)}u\circ \Phi_\iota \ dP- \int_{\mathcal C(X)}u \circ \Phi_{f_i} dP \Bigr| \le \e' P(\mathcal C(X)) \to 0 \quad (i \to \infty).
\end{align}
Since $d_P({\Phi_{f_i}}_*P, Q) \to 0$ as $i \to \infty$, we know that ${\Phi_{f_i}}_*P$ converges weakly to $Q$ as $i \to \infty$ (see e.g. \cite[\S 6]{B99}): 
\begin{align} \label{eq: ISO-2}
\int_{\mathcal C(X)}u \circ \Phi_{f_i} \ dP= \int_{\mathcal C(Y)}u \ d({\Phi_{f_i}}_*P)  \to  \int_{\mathcal C(Y)}u \ dQ \quad (i \to \infty).
\end{align}
By \eqref{eq: ISO-1} and \eqref{eq: ISO-2}, we have the equality \eqref{eq: isom} and we finish the proof.
\qed

\begin{rem} \normalfont 
When we take $X=Y$, by definition, we have $d_P(P, Q) \ge d_{LP}\bigl((X, P), (X, Q)\bigr).$ The relation between $d_P$ and $d_{LP}$ is as follows:
\begin{align} \label{eq: LP}
d_{LP}\bigl((X, P), (X, Q)\bigr)=\inf_{\substack {f:X\to X \\ \text{isometry}}}d_P({\Phi_f}_*P, Q ).
\end{align}
The following example shows that $d_{LP}\bigl((X, P), (X, Q)\bigr)=0$ does not imply $d_P(P,Q)=0$: Let $X=S^1$ with the metric $d$ where $d$ is the restriction of the Euclidean metric in $\R^2$. 
Let $x,y \in S^1$ with $x \neq y$. Let $c_x, c_y \in \mathcal C(S^1)$ denote the constant paths on $x$ and $y$, that is, $c_x(t)=x$ and $c_y(t)=y$ for all $t\in [0,T]$.
Let $\delta_{c_x}$ and $\delta_{c_y}$ be the Dirac measures on $c_x$ and $c_y$.
Let  $f: S^1 \to S^1$ be the rotation which rotate $x$ to $y$. Then, by \eqref{eq: LP}, we have $$d_{LP}((S^1, \delta_{c_x}), (S^1, \delta_{c_y})) \le d_{P}({\Phi_f}_*\delta_{c_x}, \delta_{c_y})=0.$$ We see, however,  that $d_P(\delta_{c_x}, \delta_{c_y})=d(x,y)\neq 0$.
See also Figure \ref{picture.1} below.
\begin{figure}[htbp]
\includegraphics[width=7cm, bb=0 0 454 354]{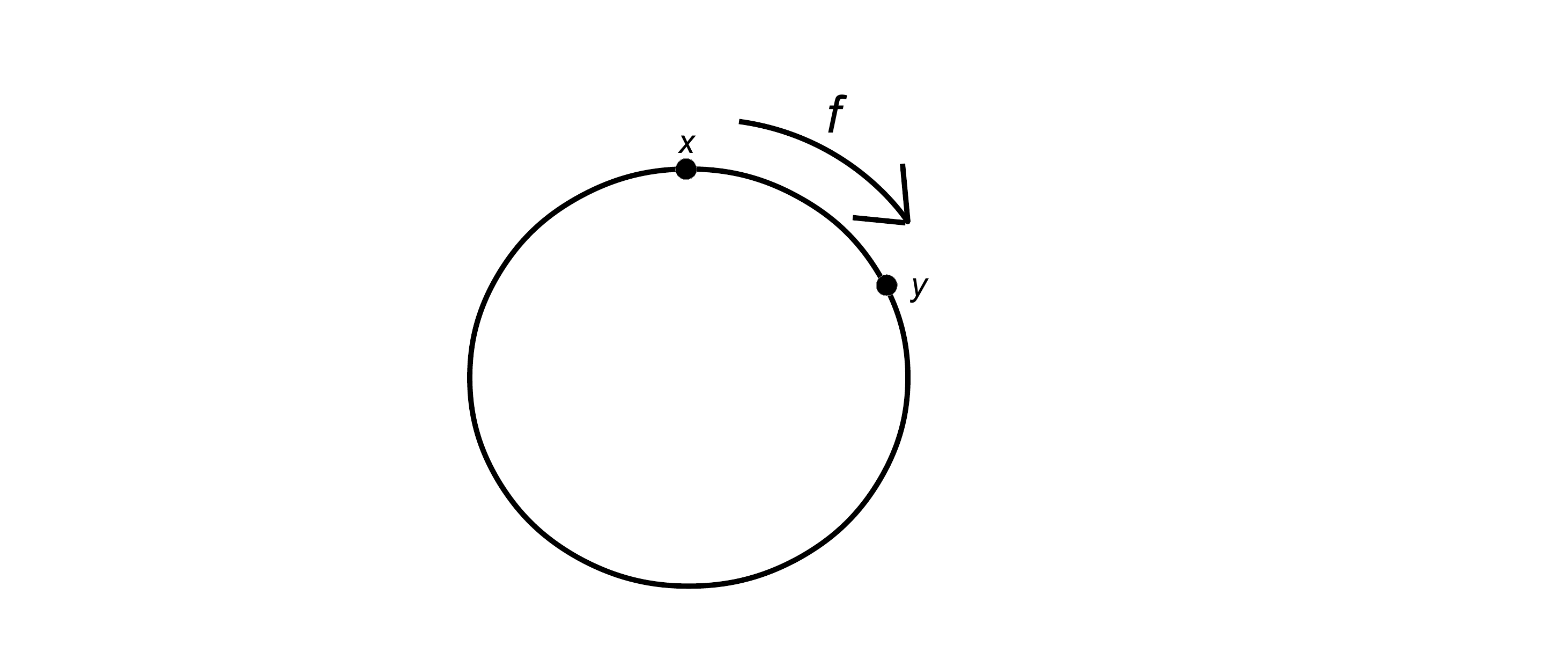}
\caption{$d_{LP}=0$ does not imply $d_P=0.$}
\label{picture.1}
\end{figure}
\end{rem}

Now we show that the metric space $(\mathcal {P}\mathcal M, d_{LP})$ is complete.
\begin{thm} \label{thm: Polish}
The metric space $(\mathcal {P}\mathcal M, d_{LP})$ is complete.
\end{thm}
\proof
Let $\{(X_i, P_i): i \in \N \}$ be a $d_{LP}$-Cauchy sequence in $\mathcal {P}\mathcal M$. It is enough to show that there are a pair $(X,P) \in \mathcal {PM}$ and a family of $(\e_i,\delta_i)$-isomorphism $f_i: (X_i,P_i) \to (X,P)$ with $\e_i,\delta_i \to 0$ as $i \to \infty$. 

{\it The existence of $X$}: The existence of $X$ follows directly from the completeness of $(\mathcal M, d_L)$. In fact, since $\{(X_i,P_i): i \in \N\}$ is a $d_{LP}$-Cauchy sequence, the sequence $\{X_i: i \in \N\}$ is a $d_L$-Cauchy sequence in $\mathcal M$. By the completeness of $(\mathcal M, d_L)$ (see Proposition \ref{prop: NSC} in Appendix), there is a compact metric space $X$ such that 
\begin{align} \label{Lip}
d_L(X_i,X) \to 0 \quad (i \to \infty).
\end{align}

{\it The existence of $P$ and $f_i$}: Since $\{(X_i, P_i): i \in \N\}$ is a $d_{LP}$-Cauchy sequence,  there is a family of $(\e_{ij}, \delta_{ij})$-isomorphisms $f_{ij}: X_i \to X_j$ for $i<j$ with $\e_{ij} \to 0$ and $\delta_{ij} \to 0$ as $i,j \to \infty$.
Take a subsequence such that $\e_{i,i+1}+\delta_{i,i+1}<1/2^{i}$.
Let ${\tilde f_{ij}}: X_i \to X_j$ be defined by 
\begin{align} \label{eq: fij}
{\tilde f_{ij}}=f_{j-1, j} \circ f_{j-2, j-1} \circ \cdot\cdot\cdot \circ f_{i,i+1} \quad (i<j),
\end{align}
and ${\tilde \e_{ij}}=\sum_{l=i}^{j-1}\e_{l,l+1},$ and ${\tilde \delta_{ij}}=\sum_{l=i}^{j-1}\delta_{l,l+1}.$ By the proof of Theorem \ref{thm: LPD}, we see that ${\tilde f_{ij}}$ is an $({\tilde \e_{ij}}, {\tilde \delta_{ij}})$-isomorphism and ${\tilde \e_{ij}}, {\tilde \delta_{ij}} \to 0$ as $i,j \to \infty$.
By the proof of Proposition \ref{prop: NSC} in Appendix (see also the equality \eqref{eq: coordinate} in Appendix), there is a family of $\e_i$-isometries $f_i:X_i \to X$ such that $\e_i\to 0$ as $i \to \infty$ and 
$$f_i\circ {\tilde f_{ji}}={f_j}.$$
This implies that  
\begin{align} \label{ineq: IJ}
\Phi_{f_i} \circ \Phi_{\tilde f_{ji}}=\Phi_{f_j}.
\end{align}

We now show that $\{{\Phi_{f_i}}_*P_i: i \in \N\}$ is a $d_P$-Cauchy sequence in $\mathcal P(\mathcal C(X))$. Let us set $$\e=\e(i,j)=({\tilde \delta_{ij}}+{\tilde \delta_{ji}}) e^{{\tilde \e_{ij}}+{\tilde \e_{ji}}+\e_i+\e_j}.$$
Note that $\e \to 0$ as $i,j \to \infty$. It suffices to show that, for any Borel set $A \subset \mathcal C(X)$, 
\begin{align*}
{\Phi_{f_i}}_*P_i (A) - {\Phi_{f_j}}_*P_j(A^{\e}) \le \e, \quad {\Phi_{f_i}}_*P_i (A^\e) - {\Phi_{f_j}}_*P_j(A) \le \e.
\end{align*}
We only show the left-hand side of the above inequalities (the right-hand side can be shown by the same argument).
For any Borel set $A \subset \mathcal C(X)$, we have  
\begin{align*}
	&{\Phi_{f_i}}_*P_i (A) - {\Phi_{f_j}}_*P_j(A^{\e})
	\\
	&= \bigl({\Phi_{f_i}}_*P_i(A)- {\Phi_{f_i \circ {\tilde f_{ji}}}}_*P_j(A^{\e})\bigr) +\bigl({\Phi_{f_i \circ {\tilde f_{ji}}}}_*P_j(A^{\e})- {\Phi_{f_j \circ {\tilde f_{ij}}}}_*P_i(A)\bigr)
	\\
	&\ + \bigl({\Phi_{f_j \circ {\tilde f_{ij}}}}_*P_i(A)- {\Phi_{f_j}}_*P_j(A^{\e}) \bigr)
	\\
	&=: {\rm(I)} +{\rm (II)}+{\rm (III)}.
\end{align*}
Since ${\tilde f_{ji}}: X_j \to X_i$ is the $({\tilde \e_{ji}}, {\tilde \delta_{ji}})$-isomorphism, we have 
\begin{align*}
{\rm(I)}  =\bigl(P_i\Phi_{f_i}^{-1}(A)- ({\Phi_{\tilde f_{ji} }}_*P_j)\Phi_{f_{i}}^{-1}(A^\e)\bigr) \le {\tilde \delta_{ji}}e^{\tilde \e_{ji}} \le \e.
\end{align*}
By the same argument, we also have ${\rm(III)}  \le \e.$
For the estimate of (II), by Lemma \ref{lem: isometry} and \eqref{ineq: IJ}, we have
\begin{align*}
	{\Phi_{f_j \circ {\tilde f_{ij}}}}_*P_i(A) & = ({\Phi_{\tilde f_{ij}}}_*P_i)\bigl(\Phi_{f_j}^{-1}(A)\bigr) \le P_j\bigl(\Phi_{f_j}^{-1}(A)^{{\tilde \delta_{ij}} e^{\tilde \e_{ij}}}\bigr) + {\tilde \delta_{ij}} e^{\tilde \e_{ij}}  
	\\
	& \le P_j\bigl(\Phi_{f_j}^{-1}(A^{{\tilde \delta_{ij}} e^{{\tilde \e_{ij}}+\e_j}}) \bigr)+ {\tilde \delta_{ij}} e^{{\tilde \e_{ij}}}
	\\
	& = P_j\bigl(\Phi^{-1}_{\tilde f_{ji}}\circ \Phi_{f_i}^{-1}(A^{\tilde \delta_{ij} e^{\tilde \e_{ij}+\e_j}}) \bigr)+ {\tilde \delta_{ij}} e^{\tilde \e_{ij}}
	\\
	&={\Phi_{f_i \circ {\tilde f_{ji}}}}_*P_j(A^{\e})+\e.
\end{align*}
We used Lemma \ref{lem: isometry} in the second line and \eqref{ineq: IJ} in the third line. 
We have ${\rm(II)}  \le \e$.
Thus $\{{\Phi_{f_i}}_*P_i: i \in \N\}$ is a $d_P$-Cauchy sequence. By the completeness of $(\mathcal P(\mathcal C(X)), d_P)$, there exists a probability measure $P$ on $\mathcal C(X)$ such that 
\begin{align}  \label{conv-1}
d_P({\Phi_{f_i}}_*P_i, P) \to 0 \quad (i \to \infty).
\end{align}

We finally show that $f_i: (X_i,P_i) \to (X,P)$ is an $(\e_i,\delta_i)$-isomorphism for some sequence $\delta_i \to 0$ as $i \to \infty$. 
By \eqref{conv-1}, there is a sequence $\delta_i \to 0$ as $i \to \infty$ such that  
\begin{align} \label{ineq: Prok-1}
{\Phi_{f_i}}_*P_i(A) \le P(A^{\delta_i})+\delta_i \quad \text{and } \quad  P(A) \le {\Phi_{f_i}}_*P_i(A^{\delta_i})+\delta_i,
\end{align}
for any Borel set $A \subset \mathcal C(X)$.
By the inequalities \eqref{ineq: Prok-1} and Lemma \ref{lem: isometry},  we have
\begin{equation} \label{ineq: Prok-2}
\begin{split}
{\Phi_{f_i^{-1}}}_*P(B) =P\bigl(\Phi^{-1}_{f_i^{-1}}(B)\bigr)
&\le {\Phi_{f_i}}_*P_i\bigl(\Phi^{-1}_{f_i^{-1}}(B)^{\delta_i}\bigr)+\delta_i 
\\
&\le P_i\bigl(\Phi^{-1}_{f_i}\circ \Phi^{-1}_{f_i^{-1}}(B^{\delta_i e^{\e_i}})\bigr) + \delta_i e^{\e_i}
\\
&\le P_i(B^{\delta_i e^{\e_i}}) + \delta_i e^{\e_i} \quad (\forall B \subset \mathcal C(X_i): \text{ Borel}).
\end{split}
\end{equation}
By the same argument, we also have 
\begin{align} \label{ineq: Prok-3}
{ P_i(B)  \le  \Phi_{f_i^{-1}}}_*P(B^{\delta_i e^{\e_i}}) + \delta_i e^{\e_i} \quad (\forall B \subset \mathcal C(X_i): \text{ Borel}).
\end{align}
Note in \eqref{ineq: Prok-1} that 
\begin{align} \label{ineq: Prok-4}
 P(A^{\delta_i})+\delta_i \le  P(A^{\delta_i e^{\e_i}})+\delta_i e^{\e_i} \quad \text{and } \quad  {\Phi_{f_i}}_*P_i(A^{\delta_i})+\delta_i \le {\Phi_{f_i}}_*P_i(A^{\delta_i e^{\e_i}})+\delta_i e^{\e_i}.
 \end{align}
Thus, by \eqref{ineq: Prok-1}, \eqref{ineq: Prok-2}, \eqref{ineq: Prok-3} and \eqref{ineq: Prok-4}, we have 
\begin{equation} \label{ineq: Prok-5}
	\begin{split}
& {\Phi_{f_i}}_*P_i(A) \le P(A^{\delta_i e^{\e_i}})+\delta_i e^{\e_i}, \ P(A) \le {\Phi_{f_i}}_*P_i(A^{\delta_i e^{\e_i}})+\delta_i e^{\e_i},
\\
&{\Phi_{f_i^{-1}}}_*P(B) \le P_i(B^{\delta_i e^{\e_i}})+\delta_i e^{\e_i}, \ P_i(B) \le {\Phi_{f_i^{-1}}}_*P_i(B^{\delta_i e^{\e_i}})+\delta_i e^{\e_i},
\end{split}
	\end{equation}
for any Borel sets $A \subset \mathcal C(X)$ and $B \subset \mathcal C(X_i)$.
Since $f_i: X_i \to X$ is an $\e_i$-isometry with $\e_i \to 0$ as $i \to 0$, the inequalities \eqref{ineq: Prok-5} means that $f_i:(X_i,P_i) \to (X,P)$ is an $(\e_i,\delta_i)$-isomorphism  with $\e_i, \delta_i \to 0$ as $i \to \infty$. We therefore have the desired result.
\qed

\begin{rem} \normalfont \label{rem: Sep}
The metric space $(\mathcal P\mathcal M, d_{LP})$ is not separable. This is because $(\mathcal M, d_L)$ is not separable as in Remark \ref{rem: NS}.
Let $X \in \mathcal M$ and $\mathcal M_X=\{Y \in \mathcal M: d_L(X,Y)<\infty\}.$ Let $\mathcal P\mathcal M_X$ be the set of isomorphism classes of pairs $(Y,P)$ where $Y \in \mathcal M_X$ and 
$P \in \mathcal P(\mathcal C(Y)).$ By \cite{Y14}, there is a $X \in \mathcal M$ such that $(\mathcal M_X, d_L)$ is not separable. Thus, for such $X$,  $(\mathcal P\mathcal M_X, d_{LP})$ is not separable.
\end{rem}

By the proof of Theorem \ref{thm: Polish}, we have the following:
\begin{cor} \label{cor: LP}
Let $(X_i, P_i), (X,P) \in \mathcal {PM}$ for all $i \in \N$.
The sequence $(X_i, P_i)$ converges to $(X,P)$ in $d_{LP}$ as $i \to \infty$ if and only if there is a family of $\e_i$-isometries $f_i:X_i \to X$ with $\e_i \to 0$ as $i \to \infty$ such that $${\Phi_{f_i}}_*P_i \to P \quad \text{weakly} \quad (i \to \infty).$$
\end{cor} 


\section{Relative compactness} \label{sec: RC}
In this section, we first give a sufficient condition for subsets in $\mathcal P\mathcal M$ to be relative compact in the case of Markov processes on Riemannian manifolds. 
Roughly speaking, this sufficient condition consists of two conditions: one is a boundedness condition of Riemannian manifolds for sectional curvatures, diameters and volumes with a fixed dimension, and the other is an upper heat kernel bound for Markov processes uniformly in each manifolds.
Second we give a sufficient condition for sequences in a relatively compact set to be convergent. This sufficient condition will be stated in terms of the generalized Mosco-convergence introduced by \cite{KS03}.

We start at  boundedness conditions for Riemannian manifolds. 
For a positive integer $n$ and positive constants $K,V,D>0$, let $\mathcal R(n, K, V, D)$ denote the set of isometry classes of $n$-dimensional connected compact Riemannian manifolds $M$ satisfying 
\begin{align} \label{condition: RC}
|{\rm sec}(M)| \le K, \quad {\rm Vol}(M) \ge V, \quad \text{and} \quad {\rm Diam}(M) \le D. 	
\end{align} 
Here ${\rm sec}(M), {\rm Vol}(M)$ and ${\rm Diam}(M)$ denote the sectional curvature, the Riemannian volume and the diameter of $M$, respectively.
We write $\mathcal R$ shortly for $\mathcal R(n, K, V, D)$.
Note that, by the Bishop inequality,  there is a constant $V'$ such that ${\rm Vol}(M)<V'$ for all $M \in \mathcal R$.
By \cite[Theorem 8.19]{Gro99} and \cite{Kat85} (see also \cite[Theorem 383]{Ber07}), we know that
\begin{align} \label{fact: LRC}
\mathcal R\text{ is relatively compact in $(\mathcal M, d_L).$}
\end{align}

We explain Markov processes considered in this section. Let $(\mathcal E, \mathcal F)$ be a regular symmetric Dirichlet form on $L^2(M, {\rm Vol})$ and $\{T_t\}_{t \in (0,\infty)}$ be the corresponding semigroup on $L^2(M, {\rm Vol})$.
We say that  $\{p(t,x,y): x,y \in M, t \in (0,\infty)\}$ is {\it a heat kernel} of $\{T_t\}_{t \in (0,\infty)}$ if $p(t,x,y)$ becomes an integral kernel of $T_t$: for $f \in L^2(M; {\rm Vol})$
$$T_tf(x)= \int_{M} f(y)p(t,x,y){\rm Vol}(dy) \quad (\forall x \in M).$$
\begin{asmp} \label{asmp: F} \normalfont
The following conditions hold:
\begin{description}
	\item[(i)]$(\mathcal E, \mathcal F)$ is a strongly local regular symmetric Dirichlet form on $L^2(M, {\rm Vol})$; 
	\item[(ii)] the corresponding semigroup $\{T_t\}_{t \in (0,\infty)}$ is a Feller semigroup and has a jointly continuous heat kernel $p(t,x,y)$  on $(0,T]\times M \times M$. 
\end{description}
\end{asmp}
See, e.g., \cite[Theorem I.9.4]{BG68} for the Feller property.
\begin{rem} \normalfont \label{rem: Feller}
We assumed the joint-continuity and the Feller property of the heat kernel only for simplicity. The following arguments can be modified by excluding null-capacity sets with respect to $(\mathcal E, \mathcal F)$
if we need to remove the Feller property assumption. The joint-continuity of a function $\phi$ which will be taken in \eqref{eq: tight} is assumed for the same reason.
\end{rem}

By \cite[Theorem I.9.4]{BG68}, there is a Hunt process 
\begin{align} \label{Hunt1}
\bigl(\Omega, \mathcal M, \{\mathcal M(t)\}_{t \in [0,\infty)}, \{{P^x}\}_{x \in M}, \{{X(t)}\}_{t\in [0,\infty)}\bigr)
\end{align}
such that, for any bounded Borel function $f$ in $L^2(M, {\rm Vol})$, we have 
$$E^x(f(X(t)))=T_tf(x) \quad (\forall t \in (0,\infty), \forall x \in M).$$
By the locality of $(\mathcal E, \mathcal F)$, we know that ${X}(\cdot)$ has continuous paths almost surely. By the strong locality and the compactness of $M$, we know that ${X}(t) \in M$ for all $t \in [0,\infty)$ almost surely. Thus we see that $X: \Omega \to \mathcal C(M)$ almost surely and the law of $X$ lives on $\mathcal C(M)$.
We refer the reader to e.g., \cite{FOT11} for details of Dirichlet forms and Hunt processes.

Let $\mu$ be a probability measure on $M$. Let $P^\mu$ denote the probability measure with the initial distribution $\mu$:
\begin{align} \label{initial}
P^{\mu}(A)=\int_AP^x\mu(dx) \quad (\forall A\subset M:\text{ Borel}).
\end{align}

Now we introduce a main object in this section, a subset $\mathcal P_\phi\mathcal R(n, K, V, D)$ of $\mathcal {PM}$ determined by a certain function $\phi$.
Let $\phi: (0,T]\times [0,D] \to [0,\infty)$ be a jointly continuous function satisfying that, for any $\e >0$,  
\begin{equation} \label{eq: tight}
\begin{split}
\lim_{\lambda \to 0}\sup_{r>\e, \xi \in (0,\lambda]}\phi(\xi, r)=0,
\end{split}
\end{equation}
where  $D>0$ is the uniform bound of diameters of elements in $\mathcal R=\mathcal R(n,K,V,D)$. 

\begin{defn} \normalfont \label{asmp: RC}
For a function $\phi$ satisfying the above conditions, 
the set $\mathcal P_\phi\mathcal R(n, K, V, D)$ is defined to be the set of isomorphism classes of  pairs $(M, P)$ where $M \in  \mathcal R$
and $P$ is  the law of $P^{\mu}$ for an initial distribution $\mu$ and a Hunt process on $M$ associated with $(\mathcal E, \mathcal F)$  satisfying Assumption \ref{asmp: F} and that 
the heat kernel $p(t,x,y)$ is dominated by $\phi$ in the following sense:  there exists a $\tau>0$ such that for all $t \in (0,\tau \wedge T]$, all $x,y \in M$ and all $M \in \mathcal R$,  
\begin{align} \label{condition: UHK}
p(t, x, y) \le \phi(t, d_M(x,y)).
\end{align}
We also write $\mathcal {P_{\phi}R}$ shortly for $\mathcal P_\phi\mathcal R(n, K, V, D)$. 
\end{defn}

Then we have the main theorem of this section:
\begin{thm} \label{thm: RRC}
The set $\mathcal {P}_{\phi}\mathcal R$ is relatively compact in $(\mathcal {P}\mathcal M, d_{LP})$.
\end{thm}
\begin{rem} \normalfont
Let $\overline{\mathcal {P}_{\phi}\mathcal R}$ be the completion of $\mathcal {P}_{\phi}\mathcal R$ with respect to $d_{LP}$.
As a byproduct of Theorem \ref{thm: RRC}, we see $$(\overline{\mathcal {P}_{\phi}\mathcal R}, d_{LP}) \text{ is a compact metric space}.$$
Let $\overline{\mathcal R }$ be the completion of ${\mathcal R }$ with respect to $d_L$.
 In general, $M \in \overline{\mathcal R }$ has a $C^{1,\alpha}$-Riemannian structure for any $0<\alpha<1$. See, e.g., \cite[Theorem 384]{Ber07}.
\end{rem}

{\it Proof of Theorem \ref{thm: RRC}:}
Since any metric spaces satisfy the first axiom of countability, it is enough to show that any sequence $(M_i, P_i)\in \mathcal {P}_{\phi}\mathcal R$ has a subsequence converging to some $(M, P) \in \mathcal {P}\mathcal M$. 
Since  $\mathcal R$ is relatively compact with respect to $d_L$ (see \eqref{fact: LRC}), the sequence $M_i \in \mathcal R$ has a converging subsequence (write also $M_i$) to a compact 
metric space $M$ with respect to $d_L$. Thus there is a family of $\e_i$-isometries $f_i: M_i \to M$ with $\e_i \to 0$ as $i \to \infty$.

By Corollary \ref{cor: LP}, the proof is completed if we show 
$\{{\Phi_{f_i}}_*P_i: i \in\N\}$ is relatively compact in $(\mathcal P(\mathcal C(M)), d_P)$.
By \cite[\S 5]{B99}, it is equivalent to show that $\{{\Phi_{f_i}}_*P_i: i \in\N\}$ is tight.
That is, for any $\e>0$, there is a compact set $K \subset \mathcal C(M)$ such that 
$${\Phi_{f_i}}_*P_i(K) >1-\e \quad (\forall i \in \N).$$
Let $(\mathcal E_i, \mathcal F_i)$ and $\mu_i$ be a strongly local regular Dirichlet form and its initial distribution associated with $P_i$.
Let a Hunt process induced by $(\mathcal E_i, \mathcal F_i)$ be denoted by 
$$(\{\mathcal M_i(t)\}_{t \in [0,\infty)}, \{P_i^x\}_{x \in M_i}, \{X_i(t)\}_{t\in [0,\infty)}).$$
By \cite[Theorem 1]{Ald78} (see also \cite[Proposition 4.3]{CKP12}), it suffices for tightness of $\{{\Phi_{f_i}}_*P_i: i \in\N\}$ to show the following two statements: 
\begin{description}
	\item[(i)] for some fixed $m \in M$, it holds that, for each positive $\eta,$ there are $a\ge 0$ and $i_0 \in \N$ such that 
	$${\Phi_{f_i}}_*P_i\bigl(v: d(v(0),m) \ge a\bigr) \le \eta \quad (\forall i \ge i_0);$$
	\item[(ii)] for any positive constants $\gamma$ and $\zeta$, there is $\lambda>0$ such that
\begin{align} \label{ineq: Ald}
	\limsup_{i \to \infty}{\Phi_{f_i}}_*P_i\bigl(v: \sup_{s: t\le s\le t+\lambda}d(v(s), v(t))>\gamma \bigr) <\zeta \quad (\forall t \in [0,T-\lambda]).
\end{align}
\end{description}

We show the statement (i). Since $f_i$ is an $\e_i$-isometry, we have that, 
for any $a>0,$
	\begin{align*}
	{\Phi_{f_i}}_*P_i\bigl(v: d(v(0),m) \ge a\bigr) & \le P_i\bigl(v: d_i(v(0),f_i^{-1}(m)) \ge ae^{-\e_i}\bigr).
	\end{align*}
Since Diam$(M)<D$ for any $M \in \mathcal R$, if we take a sufficiently large $a$ such that $\sup_iae^{-\e_i} >D$, we have $$ P_i\bigl(v: d_i(v(0),f_i^{-1}(m)) \ge ae^{-\e_i}\bigr)=0.$$
This concludes (i).

We show the statement (ii). Since ${\Phi_{f_i}}$ is an $\e_i$-isometry, we have 
\begin{align}
{\Phi_{f_i}}_*P_i\bigl(v: \sup_{s: t\le s\le t+\lambda}d(v(s), v(t))>\gamma \bigr) \le P_i\bigl(v: \sup_{s: t\le s\le t+\lambda}d_i(v(s),  v(t)) >\gamma e^{-\e_i} \bigr).
\end{align}
By the Markov property, we have, for any $t \in [0,T-\lambda]$,
\begin{equation} \label{eq: Markov1}
\begin{split}
&P^x_i\bigl(\sup_{s: t\le s\le t+\lambda}d_i(X_i(s),  X_i(t)) >\gamma e^{-\e_i} \bigr)
	\\
&=E^x_i\Bigl(P^{X_i(t)}_i\bigl(\sup_{s: t\le s\le t+\lambda}d_i(X_i(s-t),  X_i(0)) >\gamma e^{-\e_i} \bigr)\Bigr)
	\\
	& \le \sup_{x \in M_i}P^x_i\bigl(\sup_{s: t\le s\le t+\lambda}d_i(X_i(s-t),  x) >\gamma e^{-\e_i} \bigr).
\end{split}
\end{equation}
Define a $\{\mathcal M_i(t):t \in [0,T] \}$-stopping time $S_i$ as 
$$S_i=
1 \wedge \inf\{t\in [0,1): d_i(X_i(t),X_i(0))>\gamma e^{-\e_i}\}.
$$
By the strong Markov property, we have 
\begin{equation} \label{eq: SM}
\begin{split}
&P_i^x\Bigl(d_i\bigl(X_i(\lambda),X_i(S_i)\bigr)>\frac{\gamma e^{-\e_i}}{2}, S_i\le \lambda \Bigr)
\\
& =E_i^x\Bigl(\1_{S_i\le \lambda}P_i^{X_i(S_i)}\Bigl(d_i\bigl(X_i(\lambda-s),X_i(0)\bigr)>\frac{\gamma e^{-\e_i}}{2}\Bigr)\Big|_{s=S_i}\Bigr).
\end{split}
\end{equation}
By \eqref{eq: SM}, we have, for any $t \in [0,T-\lambda],$
\begin{equation} \label{ineq: Markov2}
\begin{split}
&P^x_i\Bigl(\sup_{s: t\le s\le t+\lambda}d_i\bigl(X_i(s-t),  x\bigr) >\gamma e^{-\e_i} \Bigr) = P_i^x\Bigl(S_i \le \lambda\Bigr)
\\
&\le P_i^x\Bigl(d_i\bigl(X_i(\lambda),X_i(0)\bigr)\ge\frac{\gamma e^{-\e_i}}{2}\Bigr)+P_i^x\Bigl(d_i\bigl(X_i(\lambda),X_i(0)\bigr)<\frac{\gamma e^{-\e_i}}{2}, S_i\le \lambda \Bigr)
\\
&\le P_i^x\Bigl(d_i\bigl(X_i(\lambda),X_i(0)\bigr)\ge \frac{\gamma e^{-\e_i}}{2}\Bigr)+P_i^x\Bigl(d_i\bigl(X_i(\lambda),X_i(S_i)\bigr)>\frac{\gamma e^{-\e_i}}{2}, S_i\le \lambda \Bigr)
\\
&\le P_i^x\Bigl(d_i\bigl(X_i(\lambda),X_i(0)\bigr)\ge\frac{\gamma e^{-\e_i}}{2}\Bigr)+E_i^x\Bigl(\1_{S_i\le \lambda}P_i^{X_i(S_i)}\Bigl(d_i\bigl(X_i(\lambda-s),X_i(0)\bigr)>\frac{\gamma e^{-\e_i}}{2}\Bigr)\Big|_{s=S_i}\Bigr)
\\
&\le 2 \sup_{x \in M_i, \ \xi \in (0,\lambda]}P_i^x\Bigl(d_i\bigl(X_i(\xi),X_i(0)\bigr)\ge\frac{\gamma e^{-\e_i}}{2}\Bigr)
\\
& \le 2 \sup_{x \in M_i, \ \xi \in (0,\lambda]}\int_{B(x, \gamma e^{-\e_i}/2)^c}p_i(\xi,x,y) {\rm Vol}(dy).
\end{split}
\end{equation}
Since the Riemannian volumes of $M \in \mathcal R$ are bounded by $V'$, by \eqref{eq: Markov1}, \eqref{ineq: Markov2} and \eqref{condition: UHK} of Definition \ref{asmp: RC}, for any $t\in [0,\tau\wedge T-\lambda]$ taking $\lambda$ sufficiently small as $0 <\lambda < \tau\wedge T$, we have
\begin{align*}
&P_i\Bigl(v: \sup_{s: t\le s\le t+\lambda}d_i\bigl(v(s),  v(t)\bigr) >\gamma e^{-\e_i} \Bigr)
\\
&= \int_{M_i}P^x_i\Bigl(\sup_{s: t\le s\le t+\lambda}d_i\bigl(X_i(s),  X_i(t)\bigr) >\gamma e^{-\e_i} \Bigr) \mu_i(dx)
\\
	& \le 2\sup_{x \in M_i, \ \xi \in (0,\lambda]}  \int_{B(x, \gamma e^{-\e_i}/2)^c}p_i(\xi, x, y) {\rm Vol}_{M_i}(dy)
	\\
	& \le 2\sup_{x \in M_i, \ \xi \in (0,\lambda]}\int_{B(x, \gamma e^{-\e_i}/2)^c}\phi(\xi, d_i(x,y)){\rm Vol}_{M_i}(dy)
	\\
	& \le 2\sup_{\substack{x,y \in M_i, \ \xi \in (0,\lambda], \\ d_i(x,y) \ge \gamma e^{-\e_i}/2}}\phi(\xi, d_i(x,y)) {\rm Vol}_{M_i}(B(x, \gamma e^{-\e_i}/2)^c)
	\\
	&\le 2V'\sup_{\substack{x,y \in M_i, \ \xi \in (0,\lambda], \\ d_i(x,y) > \gamma'>0 }}\phi(\xi, d_i(x,y)) <\infty.
\end{align*} 
The constant $\gamma'$ is taken to be $\gamma'=\liminf_{i \to \infty}\gamma e^{-\e_i}/2$ in the last line. The last inequality follows from the boundedness of $\phi$ on $(0,\lambda]\times[\gamma',D]$, 
which follows from the joint-continuity of $\phi$ and the condition \eqref{eq: tight}.
By \eqref{eq: tight}, we have 
\begin{align*}
&\lim_{\lambda \to 0}\limsup_{i \to \infty}\sup_{\substack{x,y \in M_i, \ \xi \in (0,\lambda],\\ d_i(x,y) > \gamma'>0}}\phi(\xi, d_i(x,y))
\\
&\le  \lim_{ \lambda \to 0}\sup_{\xi \in (0,\lambda], \ r>\gamma'>0}\phi(\xi,r)=0.
\end{align*}
This implies \eqref{ineq: Ald} and we complete the proof.
\qed

We start to consider the second objective in this section, that is, a sufficient condition for sequences in $\mathcal {P}_{\phi}\mathcal R$ to be convergent. Let $(M_i, P_i) \in \mathcal {P}_{\phi}\mathcal R$ . By Theorem \ref{thm: RRC}, we know that there is a subsequence $(M_{i'}, P_{i'})$ converging to some $(M, P)$ in the completion $\overline{\mathcal {P}_{\phi}\mathcal R}^{d_{LP}}$ with respect to $d_{LP}$.
Hereafter we consider under what conditions, the whole sequence $(M_i, P_i)$ converges to $(M, P)$. 

Let $(M_i,g_i)$ be a sequence of Riemannian manifolds with Riemannian metrics $g_i$. Assume that $M_i$ converges to some $M \in \mathcal M$ in the Lipschitz distance $d_L$ with $\e_i$-isometries $f_i:M_i \to M$. We know that the limit space $M$ has a structure of the $n$-dimensional $C^{1,\alpha}$-Riemannian manifold for any $0<\alpha<1$. That is, $M$ is a $n$-dimensional $C^{\infty}$-manifold with a $C^{1,\alpha}$-Riemannian metric $g$. See, e.g., \cite[Theorem 384]{Ber07}. Let ${\rm Vol}_i$ and ${\rm Vol}$ be Riemannian volumes induced by $g_i$ and $g$.

Let $(\E_i, \F_i)$ be a sequence of Dirichlet forms on $L^2(M_i;{\rm Vol}_i)$ and $(\E,\F)$ is a Dirichlet form on $L^2(M;{\rm Vol})$ both satisfying Assumption \ref{asmp: F}. We consider a convergence of Dirichlet forms $(\E_i,\F_i )$ to $(\E, \F)$, which is a special case of \cite[Definition 2.11]{KS03}.
Let us set $\E_i(u):=\E_i(u,u)$ for $u \in \F_i$ and $\E_i(u):=\infty$ if $u \in L^2(M_i; {\rm Vol}) \setminus \F_i$ (we also treat  $(\E, \F)$ in the same manner).
For $u \in L^2(M_i)$, we define {\it the push-forward} ${f_i}_*u \in L^2(M)$ by ${f_i}_*u(x) = u\circ f_i^{-1}(x)$ for $x \in M$.
We can check ${f_i}_*u \in L^2(M)$ because of the following inequality: 
\begin{align} \label{ineq: vol}
e^{-n\e_i}{\rm Vol} \le {f_i}_*{\rm Vol}_i \le e^{n\e_i}{\rm Vol}.
\end{align}
The above inequality follows from the definition of an $\e_i$-isometry. By this inequality, we have that  ${f_i}_*{\rm Vol}_i$ are absolutely continuous with respect to ${\rm Vol}$. 
Similarly, for $u \in L^2(M)$, we define {\it the pull-back} ${f_i}^*u$ by $u\circ f_i(x)$ for any $x \in M_i$. We can also check ${f_i}^*u \in L^2(M_i)$ by the same argument of ${f_i}_*u \in L^2(M)$.

Now we define a convergence of Dirichlet forms, which is a special case of \cite[Definition 2.11]{KS03}(see also \cite[Definition 8.1]{CKK13}).
\begin{defn} \label{defn: Mosco} \normalfont
We say that $(\E_i, \F_i)$ {\it converges in the Mosco sense} to $(\E, \F)$ if the following statement holds: there is a family of $\e_i$-isometries $f_i: M_i \to M$ with $\e_i \to 0$ as $i \to \infty$ satisfying 
\begin{description}
	\item[(i)] for any $u_i \in L^2(M_i)$ and $u \in L^2(M)$ satisfying ${f_i}_*u_i$ converges weakly to $u$ in $L^2(M)$, we have 
	$$\liminf_{i \to \infty}\E_i(u_i) \ge \E(u);$$
	\item[(ii)] for any $u \in L^2(M)$, there exists a sequence $u_i \in L^2(M_i)$ satisfying that ${f_i}_*u_i$ converges to $u$ in $L^2(M)$ and 
	$$\limsup_{i \to \infty}\E_i(u_i) \le \E(u).$$
\end{description}
\end{defn}
Note that the notion of the Mosco-convergence does not depend on a specific family of $\e_i$-isometries $f_i$ in the following sense: if $(\E_i, \F_i)$ converges in the Mosco sense to another Dirichlet form $(\E', \F')$ with respect to another family of $\e_i$-isometries $g_i: M_i \to M'$, then there is an isometry $\iota: M \to M'$ satisfying 
$$\E'(u,v)=\E(\iota^*u,\iota^*v) \quad (\forall u,v \in \D(\E')).$$

Let $\{G_i(\alpha)\}_{\alpha>0}$ and $\{G(\alpha)\}_{\alpha>0}$ be the resolvents corresponding to $(\E_i, \F_i)$ and $(\E, \F)$, respectively. We have the following statement, which is a special case 
of \cite[Theorem 2.4]{KS03} (see also \cite[Theorem 8.3]{CKK13}):
\begin{prop} \label{prop: Mosco}
The following statements are equivalent:
\begin{description}
	\item[(i)]  $(\E_i, \F_i)$ converges in the Mosco sense to $(\E, \F)$;
	\item[(ii)] there is a family of $\e_i$-isometries $f_i: M_i \to M$ satisfying $${f_i}_*T_i(t){f_i}^*u \to T(t)u \quad \text{in} \quad L^2(M),$$ for any $u \in L^2(M)$ and the convergence is uniformly in $t \in [0,T]$.
	\item[(iii)] there is a family of $\e_i$-isometries $f_i: M_i \to M$ satisfying $${f_i}_*G_i(\alpha){f_i}^*u \to G(\alpha)u \quad \text{in} \quad L^2(M),$$ for any $\alpha>0$ and any $u \in L^2(M)$.
\end{description}
\end{prop}
\proof
Modify the proof of  \cite[Theorem 8.3]{CKK13} as $\pi_iu=f_i^*u$ for $u \in L^2(M)$ and $E_iu={f_i}_*u$ for $u \in L^2(M_i)$. 
Noting the inequality \eqref{ineq: vol}, we can do the same argument in \cite[Theorem 8.3]{CKK13} and have the desired result.
\qed

\begin{asmp} \label{asmp: SC} \normalfont
For $i \in \N$, let $(M_i, P_i) \in \mathcal {P}_{\phi}\mathcal R$ for $(\mathcal E_i, \mathcal F_i)$ with an initial distribution $\mu_i=\varphi_i {\rm Vol}_i$ where $\varphi_i \in L^2(M_i).$
Let $(M, P)$ be an element in the completion $\overline{\mathcal {P_\phi R}}$ with respect to $d_{LP}$ and assume $P$ to be the law of $P^{\mu}$ for a Hunt process associated with $(\mathcal E, \mathcal F)$ satisfying Assumption \ref{asmp: F} with an initial distribution $\mu=\varphi {\rm Vol}$ with $\varphi \in L^2(M)$.
Assume that there is a family of maps $f_i: M_i \to M$ satisfying the following conditions:
\begin{description}
\item[(i)] $f_i$ is an $\e_i$-isometry with $\e_i \to 0$ as $i \to \infty$;

\item[(ii)] $f_i$ satisfies (i) and (ii) of Definition \ref{defn: Mosco};

\item[(iii)] ${f_i}_*\varphi_i$ converges to $\varphi$ in $L^2(M)$. 
\end{description}
\end{asmp}

\begin{thm} \label{thm: SCF}
If Assumption \ref{asmp: SC} holds, then $ (M_i, P_i)$ converges to $(M, P)$ as $i \to \infty$ in the sense of $d_{LP}$.
\end{thm}
{\it Proof of Theorem \ref{thm: SCF}.}
By Corollary \ref{cor: LP} and Theorem \ref{thm: RRC}, it suffices to show that the finite-dimensional distributions of ${\Phi_{f_i}}_*P_i$ converge weakly to those of $P$. 

Since $M$ is compact, any bounded continuous functions on $M$ are square-integrable, it suffices to show that, for any $k \in \N$, any $0=t_0<t_1<t_2<...<t_k\le T$ and any bounded Borel measurable functions 
$g_1, g_2, ..., g_k$ in $L^2(M)$, 
\begin{equation} \label{eq: FDC}
	\begin{split}
& E_i\bigl(f_i^*g_1\circ X_i({t_1})f_i^*g_2\circ X_i(t_2)\cdot\cdot\cdot f_i^*g_k \circ X_i(t_k)\bigr)
\\
& \to E\bigl(g_1\circ X({t_1})g_2\circ X(t_2)\cdot\cdot\cdot g_k \circ  X(t_k)\bigr) \quad (i \to \infty).
	\end{split}
\end{equation}
Let us set inductively
\begin{align*}
& h_i^k=g_k, \quad h_i^{k-1}(\cdot)=g_{k-1}(\cdot){f_i}_*T_i({t_k-t_{k-1}})f_i^*h_i^k(\cdot), ... ,
\\
& h_i^{1}(\cdot)=g_{1}(\cdot){f_i}_*T_i({t_2-t_{1}})f_i^*h_i^2(\cdot),
\end{align*}
and
\begin{align*}
& h^k=g_k, \quad h^{k-1}(\cdot)=g_{k-1}(\cdot)T({t_k-t_{k-1}})h^k(\cdot), ... ,
\\
& h^{1}(\cdot)=g_{1}(\cdot)T({t_2-t_{1}})h^2(\cdot).
\end{align*}
By Proposition \ref{prop: Mosco} and boundedness of $g_k$, we have
\begin{align} \label{conv: FD1}
	\|h^{k-1}_i-h^{k-1}\|_{L^2(M)}&=\|  g_{k-1}(\cdot){f_i}_*T_i({t_k-t_{k-1}})f_i^*h^k_i(\cdot)- g_{k-1}(\cdot)T({t_k-t_{k-1}})h^k(\cdot)\|_{L^2(M)}  \notag
	\\
	&\to 0 \quad (i \to \infty).
\end{align}
Inductively, we have 
\begin{align*}
&\|h^{k-2}_i-h^{k-2}\|_{L^2(M)}
\\
&=\|g_{k-2}(\cdot){f_i}_*T_i(t_{k-1}-t_{k-2})f_i^*h_i^{k-1}(\cdot)-g_{k-2}(\cdot)T(t_{k-1}-t_{k-2})h^{k-1}(\cdot)\|_{L^2(M)} 
\\
& \le \| g_{k-2}{f_i}_*T_i(t_{k-1}-t_{k-2})f_i^*h_i^{k-1}-g_{k-2}{f_i}_*T_i(t_{k-1}-t_{k-2})f_i^*h^{k-1}\|_{L^2(M)}
\\
&\quad + \|g_{k-2}{f_i}_*T_i(t_{k-1}-t_{k-2})f_i^*h^{k-1}- g_{k-2}T(t_{k-1}-t_{k-2})h^{k-1}\|_{L^2(M)}
\\
& \le \|g_{k-2}\|_{\infty}\| {f_i}_*T_i(t_{k-1}-t_{k-2})\|_{\text{op}}\|f_i^*h_i^{k-1}-f_i^*h^{k-1}\|_{L^2(M)}
\\
&\quad + \|g_{k-2}\|_{\infty}\|{f_i}_*T_i(t_{k-1}-t_{k-2})f_i^*h^{k-1}- T(t_{k-1}-t_{k-2})h^{k-1}\|_{L^2(M)}
\\
&=:{\rm (I)}_i + {\rm (II)}_i,
\end{align*}
where $\|{f_i}_*T_i(t_{k-1}-t_{k-2})\|_{\text{op}}$ means the operator norm of ${f_i}_*T_i(t_{k-1}-t_{k-2}): L^2(M_i) \to L^2(M)$. 

The quantity (II)$_i$ converges to $0$ as $i \to \infty$ by (ii) of Assumption \ref{asmp: SC} and Proposition \ref{prop: Mosco}. 

We estimate (I)$_i$.
By the inequality \eqref{ineq: vol} and the contraction property of the semigroup $\{T_i(t)\}_{t>0}$, we can check easily that there is a constant $C$ independent to $i$ satisfying 
\begin{align} \label{conv: FD2}
\| {f_i}_*T_i(t_{k-1}-t_{k-2})\|_{\text{op}} \le C.
\end{align} 
By \eqref{conv: FD1} and the inequality \eqref{ineq: vol}, we have 
\begin{align} \label{conv: FD3}
\|f_i^*h_i^{k-1}-f_i^*h^{k-1}\|_{L^2(M)} \to 0 \quad (i \to \infty).
\end{align}
Thus we have ${\rm (I)}_i \to 0$ as $i \to \infty$.

By using the above argument inductively and the Markov property, we have 
\begin{align} \label{eq: FDC2}
&\|h^{1}_i-h^{1}\|_{L^2(M)} \notag
\\
& = \Bigl\|E^{f_i^{-1}(x)}_i\bigl(f_i^*g_1\circ X_i({t_1})f_i^*g_2\circ X_i(t_2)\cdot\cdot\cdot f_i^*g_k \circ X_i(t_k)\bigr) \notag
\\
& \quad - E^x\bigl(g_1\circ X({t_1})g_2\circ X(t_2)\cdot\cdot\cdot g_k \circ  X(t_k)\bigr)\Bigr\|_{L^2(M)}  \notag
\\
&\to 0 \quad (i \to \infty).
	\end{align}
On the other hand, by the inequality \eqref{ineq: vol}, we have that 
\begin{align} \label{con: uniform}
	\frac{d({f_i}_*{\rm Vol}_i)}{d {\rm Vol}} \to \mathbf{1}_M \quad \text{uniformly},
\end{align}
where $\mathbf{1}_M$ means the indicator function on $M$.
By the fact \eqref{con: uniform} and (iii) of Assumption \ref{asmp: SC}, we have
 \begin{align} \label{conv: ID}
\| \frac{d({f_i}_*(\varphi_i {\rm Vol}_i))}{d {\rm Vol}}-\varphi\|_{L^2(M)} \to 0 \quad (i \to \infty).
 \end{align}
Thus, by \eqref{eq: FDC2} and \eqref{conv: ID}, using the Schwarz inequality, we have 
\begin{align*}
	&\bigl|E_i\bigl(f_i^*g_1\circ X_i({t_1})f_i^*g_2\circ X_i(t_2)\cdot\cdot\cdot f_i^*g_k \circ X_i(t_k)\bigr)
	\\
	& \quad - E\bigl(g_1\circ X({t_1})g_2\circ X(t_2)\cdot\cdot\cdot g_k \circ  X(t_k)\bigr) \bigr| \quad (i \to \infty) 
	\\
	&= \Bigl|\int_{M}E^{f_i^{-1}(x)}_i\bigl(f_i^*g_1\circ X_i({t_1})f_i^*g_2\circ X_i(t_2)\cdot\cdot\cdot f_i^*g_k \circ X_i(t_k)\bigr) \ {f_i}_*(\varphi_i {\rm Vol}_i)(dx)
	\\
	& \quad - \int_{M}E^x\bigl(g_1\circ X({t_1})g_2\circ X(t_2)\cdot\cdot\cdot g_k \circ  X(t_k)\bigr) \varphi(x){\rm Vol}(dx) \Bigr|
	\\
	& = \Bigl|\int_{M}E^{f_i^{-1}(x)}_i\bigl(f_i^*g_1\circ X_i({t_1})f_i^*g_2\circ X_i(t_2)\cdot\cdot\cdot f_i^*g_k \circ X_i(t_k)\bigr) \ \frac{d({f_i}_*(\varphi_i {\rm Vol}_i))}{d {\rm Vol}}(x){\rm Vol}(dx)
	\\
	& \quad - \int_{M}E^x\bigl(g_1\circ X({t_1})g_2\circ X(t_2)\cdot\cdot\cdot g_k \circ  X(t_k)\bigr) \varphi(x){\rm Vol}(dx) \Bigr|
	\\
	&\to 0 \quad (i \to \infty).
\end{align*}
We therefore have shown \eqref{eq: FDC} and we have completed the proof.
\qed

\section{Examples} \label{sec: EX}
\subsection{Brownian motions on Riemannian manifolds}  \label{subsec: BM}
In this subsection, we consider the case when a state space $M$ is in $\mathcal R=\mathcal R(n,K,V,D)$ and a probability measure $P$ is the law of the Brownian motion, which is the Markov process induced by the Laplacian on $M$.
In this case, the convergence of processes should follow only from the convergence of state spaces. 
We show that the convergence in $d_{LP}$ follows only from the convergence in $d_L$.

Let $(M,g)$ be in $\mathcal R$. Let $\nabla$ denote the gradient operator induced by $g$.
Let $(\mathcal E, \mathcal F)$  be the smallest closed extension of the following bilinear form on $L^2(M; {\rm Vol})$:
\begin{align} \label{eq: SEF}
\mathcal E (u,v)=\frac{1}{2}\int_Mg_x(\nabla u,\nabla v)\frac{1}{{\rm Vol}(M)}{\rm Vol}(dx) \quad (u,v \in C^{\infty}(M)).
\end{align}
We write $\overline{\rm Vol}(dx)={\rm Vol}(dx)/{\rm Vol}(M)$.

\begin{defn} \normalfont
The set $\mathcal {LR}(n,K,V,D)$ is defined to be the set of isomorphism classes of pairs $(M, P)$ where $M \in \mathcal R$ and $P$ is the law of $P^{\mu}$ for a Markov process on a time interval $[0,T]$ associated with  $(\mathcal E, \mathcal F)$ defined in \eqref{eq: SEF} with an initial probability measure $\mu=\overline{{\rm Vol}}_M$.
We denote $\mathcal {LR}$ shortly for $\mathcal {LR}(n,K,V,D)$.
\end{defn}

We show the relative compactness of $\mathcal {LR}$.
\begin{prop} \label{prop: LR}
The set $\mathcal {LR}$ is relatively compact in $(\mathcal{PM}, d_{LP})$.
\end{prop}

Before the proof, we recall the uniform heat kernel estimate of \cite[\S 4]{KK96}.
Let $\overline{p_M}(t,x,y)$ be a heat kernel of the standard energy form $(\E, \F)$ with respect to the normalized volume measure $\overline{{\rm Vol}}$ for $M \in \mathcal R$.
We know that $\overline{p_M}(t,x,y)$ is jointly continuous in $(t,x,y)$ and has the Feller property (see, e.g.,  \cite[Theorem 7.16 \& 7.20]{Gri09}).
By \cite[\S 4]{KK96}, we have the following heat kernel estimate:
\begin{align} \label{ineq: UHKE}
\overline{p_M}(t,x,y) \le \frac{C}{t^{1+\nu}}\exp \Bigl(  \frac{d_M(x,y)^2}{4t}\Bigr),
\end{align}
for all $x,y \in M$ and $0< t \le D^2$, where $\nu=\nu(n,K,D)>2$ and $C=C(n,K,D)>0$ are positive constants depending only on $n,K$ and $D$.
The important point is that $\nu$ and $C$ do not depend on each $M \in \mathcal R$.
Note that if we have $|{\rm Sec}(M)| \le K$ for $K > 0$,  we have ${\rm Ric}(M) \ge -K$. 

Now we show Proposition \ref{prop: LR}. 

{\it Proof of Proposition \ref{prop: LR}.}
Let $p_M(t,x,y)$ be the heat kernel of the standard energy form $(\E,\F)$ with respect to the Riemannian volume measure ${\rm Vol}$ (not with respect to $\overline{\rm Vol}$).
Note that 
\begin{align} \label{eq: N=N}
\frac{\overline{p_M}(t,x,y)}{{\rm Vol}(M)}=p_M(t,x,y).
\end{align}

It suffices to show that there is a jointly continuous function $\phi :(0,T]\times [0,D] \to [0,\infty)$ satisfying \eqref{eq: tight} and dominating $p_M(t,x,y)$ as \eqref{condition: UHK}.
In fact, if we show this, we have 
$$\mathcal {LR}\subset \mathcal {P}_{\phi}\mathcal {R}.$$
Since $\mathcal {P}_{\phi}\mathcal {R}$ is relatively compact by Theorem \ref{thm: RRC}, we obtain the desired result.
The existence of $\phi$ satisfying \eqref{eq: tight} and \eqref{condition: UHK} follows from \cite{KK96}. In fact, 
by  \eqref{ineq: UHKE}, \eqref{eq: N=N} and the lower-bounds $V$ of the volumes,  we have that there is a constant $C'=C'(n,K,V,D)>0$ such that 
\begin{align} \label{ineq: HKB2}
p_M(t,x,y) \le \frac{C'}{t^{1+\nu}}\exp\Bigl(-\frac{d_M(x,y)^2}{4t}\Bigr),
\end{align} 
for all $M \in \mathcal R$, all $t \in (0,D^2]$ and all $x,y \in M$. Note that the constant $C'$ does not depend on each $M \in \mathcal R$.
Thus we have checked \eqref{condition: UHK} with $\tau=D^2$.
Let $$\phi(\xi,r)=\frac{C'}{\xi^{1+\nu}}\exp(-\frac{r^2}{4\xi}).$$ Then we can check easily that $\phi$ satisfies \eqref{eq: tight}.
Thus we have completed the proof.
\qed

Let $(M_i, P_i) \in \mathcal {LR}$ and assume $M_i$ converges to some $M \in \overline{\mathcal R}^{d_L}$ where $\overline{\mathcal R}^{d_L}$ denotes the completion of $\mathcal R$ with respect to $d_L$.
As stated in Section \ref{sec: RC}, the limit space $M$ has a $C^{1,\alpha}$-Riemannian structure.
Such manifolds are in the class of Lipschitz--Riemannian manifolds (see, e.g., \cite[\S 3]{KS03}). In this framework, we have a Riemannian volume ${\rm Vol}_M$ induced by a $C^{1,\alpha}$-Riemannian metric $g$
and the standard energy form $(\mathcal E, \F)$ defined by a similar way to \eqref{eq: SEF} with respect to the weak derivative. See the detail in \cite[\S 3]{KS03} and references therein.

Let $P$ be a law of Markov process on $M$ associated with the above $(\mathcal E, \mathcal F)$ whose initial distribution is the Riemannian volume $\overline{{\rm Vol}_M}$.
Then we have the following:
\begin{prop} \label{prop: DLP}
If $M_i$ converges to $M$ in $d_L$, then $(M_i, P_i)$ converges to $(M,P) $ in $d_{LP}$.
\end{prop}
\proof
By Theorem \ref{thm: SCF}, it is sufficient to check that Assumption \ref{asmp: SC} are satisfied.
Since we assume that $M_i \to M$ in $d_{L}$, there is a family of $\e_i$-isometries  $f_i: M_i \to M$ with $\e_i \to 0$ as $i \to \infty$.
By the inequality \eqref{ineq: vol}, we can easily check that there is a function $h: [0,\infty)\to [0,\infty)$ with $\lim_{r \to 0}h(r)=0$ satisfying the following inequalities (see \cite[\S 3]{KS03}): for any $u \in L^2(M)$ and $u_i \in L^2(M_i)$, 
\begin{align*} 
&\bigl|\|f_i^*u\|_{L^2(M_i)}-\|u\|_{L^2(M)}\bigr| \le h(\e_i)\|u\|_{L^2(M)},
\\
&\bigl|\E_i(f_i^*u)-\E(u)\bigr| \le h(\e_i)\E(u), 
\\
&\bigl|\E_i(u_i)-\E({f_i}_*u_i)\bigr| \le h(\e_i)\E_i(u_i).
\end{align*}
By these inequalities, we can check easily that the conditions of Definition \ref{defn: Mosco} are satisfied with $f_i$ (see \cite[Proposition 3.1]{KS03}). 
Since we take $\mu_i={\rm Vol}_{M_i}$ and $\mu={\rm Vol}_M$ in this section, there is nothing to check about (iii) of Assumption \ref{asmp: SC}.
Thus the conditions in Assumption \ref{asmp: SC} are satisfied and we finish the proof.
\qed

By Proposition \ref{prop: DLP}, we know what is the completion $\overline{\mathcal{LR}}^{d_{LP}}$ of $\mathcal{LR}$ with respect to $d_{LP}$.
We define the subset $\mathcal {L}\overline{\mathcal{R}}^{d_L}$ consisting of pairs $(M,P)$ where $M \in \overline{\mathcal R}^{d_L}$ and $P$ is the law of $P^{\mu}$ where 
$P$ is the Brownian motion associated with the standard energy form $(\mathcal E, \mathcal F)$ defined in \eqref{eq: SEF} with the initial distribution $\mu=\overline{{\rm Vol}}_M$.

\begin{cor} We have the following:
\begin{align*}
	\overline{\mathcal{LR}}^{d_{LP}}=\mathcal {L}\overline{\mathcal{R}}^{d_L}.
\end{align*}
\end{cor}
\proof
We first show $\overline{\mathcal{LR}}^{d_{LP}}\subset \mathcal {L}\overline{\mathcal{R}}^{d_L}.$ Let $(M,P) \in \overline{\mathcal{LR}}^{d_{LP}}$. Then we have a sequence $(M_i,P_i)\in \mathcal {LR}$ such that 
 $(M_i,P_i) \to (M,P)$ in $d_{LP}$. Thus $M_i\in \mathcal {R}$ converges to $M$ in $d_{L}$, and $M \in \overline{\mathcal R}^{d_L}$.  
By Proposition \ref{prop: DLP}, we have that $P$ is the law of the Brownian motion associated with the standard energy form on $M$ with the initial distribution $\mu=\overline{{\rm Vol}}_M$.
Thus $(M,P) \in \mathcal {L}\overline{\mathcal{R}}^{d_L}$.

We second show $\overline{\mathcal{LR}}^{d_{LP}}\supset \mathcal {L}\overline{\mathcal{R}}^{d_L}.$ Let $(M,P) \in\mathcal {L}\overline{\mathcal{R}}^{d_L}.$
Then we have a sequence $M_i \to M$ in $d_L$. Let $P_i$ be the law of the Brownian motion on $M_i$ associated with the standard energy form with the initial distribution 
 $\mu=\overline{{\rm Vol}}_M$. Thus, by Proposition \ref{prop: DLP}, we have that $(M_i, P_i) \to (M,P)$ in $d_{LP}$ and thus $(M,P) \in \overline{\mathcal{LR}}^{d_{LP}}$. 
\qed

\subsection{Uniformly elliptic diffusions on Riemannian manifolds} \label{subsec: UE}
In this subsection, we consider $(M,P)$ where $(M,g) \in \mathcal R$ and $P$ is a law of a Markov process associated with another smooth Riemannian metric $h$ comparable to the given Riemannian metric $g$, that is, 
there is a constant $\Lambda>1$ satisfying 
\begin{align} \label{ineq: um}
\Lambda^{-1}g \le h \le \Lambda g. 
\end{align}
The generator associated with $h$ is a second order differential operator having smooth coefficients with the uniform elliptic condition in local coordinates.

To be precise, let $\nabla_h$ denote the gradient operator induced by $h$ satisfying \eqref{ineq: um}.
Let ${\rm Vol}_h$ be the volume measure associated with $h$.
Let $(\mathcal E^h, \mathcal F^h)$  be the smallest closed extension of the following bilinear form on $L^2(M; {\rm Vol}_h)$:
\begin{align} \label{eq: SEF2}
\mathcal E^h (u,v)=\frac{1}{2}\int_Mh_x(\nabla_h u,\nabla_h v)\frac{1}{{\rm Vol}_h(M)}{\rm Vol}_h(dx) \quad (u,v \in C^{\infty}(M)).
\end{align}
We write $\overline{\rm Vol}_h(dx)={\rm Vol}_h(dx)/{\rm Vol}_h(M)$.

\begin{defn} \normalfont
For $\Lambda >1$,  the set $\mathcal {L}_\Lambda\mathcal {R}(n,K,V,D)$ is defined to be the set of isomorphism classes of pairs $(M, P)$ where $(M,g) \in \mathcal R$ and $P$ is the law of $P^{\mu}$ for a Markov process on $[0,T]$ associated with  $(\mathcal E^h, \mathcal F^h)$ defined in \eqref{eq: SEF2} for $h$ satisfying \eqref{ineq: um} with an initial probability measure $\mu=\overline{{\rm Vol}}_h$.
We denote $\mathcal {L}_\Lambda\mathcal {R}$ shortly for $\mathcal {L}_\Lambda \mathcal {R}(n,K,V,D)$.
\end{defn}

We show the relative compactness of  $\mathcal {L}_\Lambda\mathcal {R}$.
\begin{prop} \label{prop: LR1}
The set  $\mathcal {L}_\Lambda\mathcal {R}$ is relatively compact in $(\mathcal{PM}, d_{LP})$.
\end{prop}
\proof
By \cite[\S 4]{KK96}, we have the same heat kernel estimate as \eqref{ineq: UHKE} for $(\mathcal E^h, \mathcal F^h)$. Note that, of course,  the positive constant $C$ in \eqref{ineq: UHKE} depends also on $\Lambda$ in this case.
 Thus the proof follows from the same argument of Proposition \ref{prop: LR}.
\qed

\appendix \section{Appendix}
Recall that $\mathcal M$ is the set of isometry classes of compact metric spaces and $d_L$ is the Lipschitz distance on $\mathcal M$ (see Section \ref{sec: LP}).
We show the completeness of the metric space $(\mathcal M, d_L)$.
\begin{prop} \label{prop: NSC}
$(\mathcal M, d_L)$ is a complete metric space.
\end{prop}
\proof
Let $\{X_i: i \in \N\}$ be a $d_L$-Cauchy sequence in $\mathcal M$. It suffices to show that there are a compact metric space $X \in \mathcal M$ and $\e_i$-isometries $f_i:X_i \to X$ with $\e_i \to 0$ as $i \to \infty$.

{\it The construction of $X$}: Let $f_{ij}: X_i \to X_j$ be an $\e_{ij}$-isometry for $i<j$ where $\e_{ij} \to 0$ as $i,j \to \infty.$  
Take a subsequence such that $\e_{i,i+1}<1/2^{i}$.
Let ${\tilde f_{ij}}: X_i \to X_j$ be defined by 
\begin{align} \label{eq: fij}
{\tilde f_{ij}}=f_{j-1, j} \circ f_{j-2, j-1} \circ \cdot\cdot\cdot \circ f_{i,i+1} \quad (i<j),
\end{align}
and ${\tilde \e_{ij}}=\sum_{l=i}^{j-1}\e_{l,l+1}.$ Then ${\tilde f_{ij}}$ is an ${\tilde \e_{ij}}$-isometry and ${\tilde \e_{ij}} \to 0$ as $i,j \to \infty$.
Since every compact metric space is separable, there is a countable dense subset 
$\{x^1_\alpha: \alpha \in \N\} \subset X_1$.  
We define, for any $i >1$, $$x_\alpha^i={\tilde f_{1i}}(x_\alpha^1).$$
Since ${\tilde f_{1i}}$ is a homeomorphism, the subset $\{x_\alpha^i: i \in \N\}$ is dense in $X_i$ for each $i$.
Fix $\alpha, \beta \in \N$, and consider the sequence of the real numbers 
\begin{align} \label{seq: B}
\{d(x^i_\alpha,x_\beta^i) : i\in \N\}.
\end{align}
Since $\{{\tilde f_{1i}}:i \in \N\}$ has a bounded Lipschitz constant and the compact metric space $X_1$ is bounded, we have that \eqref{seq: B} is a bounded sequence:
$$d(x^i_\alpha,x_\beta^i )\le \sup_{i}{\rm dil}(\tilde{f}_{1i})d(x^1_\alpha,x^1_\beta) < \infty.$$
Thus we can take a subsequence of \eqref{seq: B} converging to some 
real number, write $r(\alpha, \beta)$.
We can check that $r$ becomes a metric on $\{\alpha: \alpha \in \N\}$. In fact, if $r(\alpha, \beta)=0$, we have $\alpha=\beta$ because 
\begin{align} \label{ineq: inf}
0<\frac{1}{\sup_i{\rm dil}({\tilde f_{1i}}^{-1})}d(x_\alpha^1, x_\beta^1) \le d(x_\alpha^i, x_\beta^i).
\end{align}
By definition, $r(\alpha, \alpha)=0$ and $r$ is symmetric and non-negative. It is easy to see the triangle inequality.
Let $(X,d)$ be the completion of the metric space $(\{\alpha: \alpha \in \N\}, r)$. The compactness of $(X,d)$ will be shown in the next paragraph.
 
{\it The construction of $\e_i$-isometries $f_i$}: We define a map $f_i: \{x_\alpha^i: \alpha \in \N\} \to X$ by 
$$f_{i}(x_\alpha^i)=\alpha.$$
Now we extend the map $f_i$ to the whole space $X_i$. Since ${\rm dil}(\tilde{f}_{ij})$ is bounded, we have
\begin{align}\label{ineq: infsup}
d(f_i(x_\alpha^i), f_i(x^i_\beta))&= d(\alpha, \beta)=\lim_{j \to \infty}d(x_\alpha^j,x_\beta^j)=\lim_{j \to \infty}d({\tilde{f}_{ij}}(x_\alpha^i),{\tilde{f}_{ij}} (x_\beta^i)) \notag
\\
&\le \sup_{j}{\rm dil}({\tilde f_{ij}})d(x_\alpha^i, x_\beta^i) < \infty.
\end{align}
Let $x_{\alpha(n)}^i \to x^i \in X_i$ as $n \to \infty$. By the inequality \eqref{ineq: infsup}, we have that $\lim_{n \to \infty}f_i(x_{\alpha(n)}^i)$ exists.
This limit does not depend on the way of taking sequences converging to $x^i$ (use the triangle inequality to check it). Thus we define $$f_i(x^i)=\lim_{n \to \infty}f_i(x_{\alpha(n)}^i),$$ and this is well-defined. 
Thus we have extended the map $f_i$ to the whole space $X_i$.

Now we check that $f_i$ is bi-Lipschitz. We have
\begin{align}\label{ineq: infsup2}
d(f_i(x_\alpha^i), f_i(x^i_\beta))&= d(\alpha, \beta)=\lim_{j \to \infty}d(x_\alpha^j,x_\beta^j)=\lim_{j \to \infty}d({\tilde{f}_{ij}}(x_\alpha^i),{\tilde{f}_{ij}} (x_\beta^i)) \notag
\\
&\ge \frac{1}{\sup_{j}{\rm dil}({\tilde f^{-1}_{ij}})}d(x_\alpha^i, x_\beta^i).
\end{align} 
Note that $0<\sup_{j}{\rm dil}({\tilde f^{-1}_{ij}})<\infty$.  By the inequality \eqref{ineq: infsup2}, we see that $f_i$ is bijective.
By the inequality \eqref{ineq: infsup} and \eqref{ineq: infsup2}, we have $f_i$ is bi-Lipschitz.
Since $f_i$ is a homeomorphism and $X_i$ is compact, we see that $X=f_i(X_i)$ is compact. Thus $X \in \mathcal M$.

Finally we check that $f_i$ is an $\e_i$-isometry for some $\e_i \to 0$ as $i \to \infty$. We set $$\e_i=\max\{|\log(\sup_{j}{\rm dil}({\tilde f_{ij}}^{-1}))|, |\log(\sup_{j}{\rm dil}({\tilde f_{ij}}))|\}.$$ 
Then, by the inequality \eqref{ineq: infsup} and \eqref{ineq: infsup2}, we can see that $f_i: X_i \to X$ is an $\e_i$-isometry with $\e_i\to 0$ as $i \to \infty$. 
Thus we have shown that $X$ is the $d_L$-limit of $X_i$.
We have completed the proof.
\qed

Note that, in the above proof, we have that 
\begin{align} \label{eq: coordinate}
f_j\circ {\tilde f_{ij}}={f_i}.
\end{align}
We use \eqref{eq: coordinate} in the proof of Theorem \ref{thm: Polish}.

\begin{rem} \normalfont \label{rem: NS}
Note that $(\mathcal M, d_L)$ is not separable. This is  because of the following two facts: 
\begin{description}
\item[(a)]if $d_L(X,Y)<\infty$, the Hausdorff dimensions of $X$ and $Y$ must coincide;
\item[(b)]for any non-negative real number $d$, there is a compact metric space $X$
whose Hausdorff dimension is equal to $d$. 
\end{description}
See, e.g., \cite[Proposition 1.7.19]{BBI01} for (a) and \cite{SS} for (b).
Let $X \in \mathcal M$ and $\mathcal M_X=\{Y \in \mathcal M: d_L(X,Y)<\infty\}.$ We also note that there is a $X \in \mathcal M$ such that even when we restrict $d_L$ to $\mathcal M_X$, the metric space $(\mathcal M_X, d_L)$ is not separable. See \cite{Y14}.
\end{rem}

\section*{Acknowledgment}
I would like to express my great appreciation to Yohei Yamazaki for a lot of valuable and constructive comments and useful discussions. 
I would also like to thank Prof.\ Kouji Yano for careful reading of my manuscript and successive encouragement.

\end{document}